\documentclass[12pt]{article}
\usepackage{amsmath,amsfonts,amssymb}
\usepackage[pagebackref=true]{hyperref}
\def\C{\mathbb{C}}
\def\a{\mathbf{a}}
\def\b{\mathbf{b}}

\def\Be{\mathcal{B}}

\def\R{\mathbb{R}}
\def\Z{\mathbb{Z}}
\def\L{\mathcal{L}}
\def\A{\mathbf{A}}
\def\B{\mathbf{B}}
\def\L{\mathcal{L}}

\newtheorem{theorem}{\hspace*{\parindent}Theorem}
\newtheorem{lemma}{\hspace*{\parindent}Lemma}
\newtheorem{corollary}{\hspace*{\parindent}Corollary}
\newtheorem{conjecture}{\hspace*{\parindent}Conjecture}
\newtheorem{open}{\hspace*{\parindent}Open problem}

\title{The Fox-Wright function near the singularity and branch cut}
\author{D.B.\:Karp$^{\rm a}$\footnote{Corresponding author. E-mail: D.B.\:Karp -- \emph{dmitriibkarp@tdtu.edu.vn}, E.G.Prilepkina --  \emph{pril-elena@yandex.ru}}~~and E.G.Prilepkina$^{\rm b,c}$
\\[10pt]
\small{\textit{$\phantom{1}^a$Ton Duc Thang University,  Ho Chi Minh City, Vietnam}}\\
\small{\textit{$\phantom{1}^b$Far Eastern Federal University,
Vladivostok, Russia}}\\\small{\textit{$\phantom{1}^c$Institute of
Applied Mathematics, FEBRAS}}}
\date{}

\begin{document}
\maketitle

\begin{abstract}
The Fox-Wright function is a further extension of the generalized hypergeometric function obtained by introducing arbitrary positive scaling factors into the arguments of the gamma functions in the summand.  Its importance comes mostly from its role in fractional calculus although other interesting applications also exist. If the sums of the scaling factors in the top and bottom parameters are equal, the series defining the Fox-Wright function has a finite non-zero radius of convergence.  It was demonstrated by Braaksma in 1964 that the Fox-Wright function can then be extended to a holomorphic function in the complex plane cut along a ray from the positive point on the boundary of the disk of convergence to the point at infinity.  In this paper we study the behavior of the Fox-Wright function in the neighborhood of this positive singular point.  Under certain restrictions we give a convergent expansion with recursively computed coefficients completely characterizing this behavior.  We further compute the jump and the average value of the Fox-Wright function on the banks of the branch cut.
\end{abstract}

\bigskip


\bigskip


\bigskip

\section{Introduction and preliminaries}
We will use the standard symbols $\Z$, $\R$ and $\C$ to denote the sets of integer, real, and complex numbers, respectively.
We will further employ the self-explanatory notation like $\Z_{\ge2}$ or $\Z_{<0}$ for subsets of $\Z$.
Given positive vectors $\A=(A_1,\ldots,A_p)$, $\B=(B_1,\ldots,B_q)$ and complex vectors $\a=(a_1,\ldots,a_p)$, $\b=(b_1,\ldots,b_q)$,
the Fox-Wright function (or ''the generalized Wright function'') is defined by the series \cite[\S12]{Braaksma}, \cite{Wright35}, \cite[Section~4.1]{HTF1},
\begin{equation}\label{eq:pPsiqdefined}
{_{p}\Psi_q}\left(\left.\!\!\begin{array}{c}(a_1,A_1),\ldots,(a_p,A_p)\\(b_1,B_1),\ldots,(b_q,B_q)\end{array}\right|z\!\right)
={_{p}\Psi_q}\left(\left.\!\!\begin{array}{c}(\a,\A)\\(\b,\B)\end{array}\right|z\!\right)
=\sum\limits_{n=0}^{\infty}\frac{\Gamma(\A n+\a)}{\Gamma(\B n+\b)}\frac{z^n}{n!},
\end{equation}
where $\Gamma$ stands for Euler's gamma function and the shorthand notation $\Gamma(\A{n}+\a)=\prod_{j=1}^{p}\Gamma(A_{j}n+a_{j})$ (similarly for $\Gamma(\B{n}+\b)$) has been used.  We adopt this abbreviation throughout the paper.  The series (\ref{eq:pPsiqdefined}) has a nonzero radius of convergence if
\begin{equation}\label{eq:balance}
\Delta:=\sum\nolimits_{j=1}^{q}B_j-\sum\nolimits_{i=1}^{p}A_i\ge-1.
\end{equation}
More precisely, if $\Delta>-1$ the series converges for all finite values of $z$ to an entire function, while for $\Delta=-1$, its radius of convergence equals
\begin{equation}\label{eq:rho}
\rho:=\prod\limits_{k=1}^{p}A_k^{-A_k}\prod\limits_{j=1}^{q}B_j^{B_j}.
\end{equation}
Convergence on the boundary $|z|=\rho$ depends on the value of
\begin{equation}\label{eq:mu}
\mu:=\sum\nolimits_{j=1}^{q}b_j-\sum\nolimits_{k=1}^{p}a_k+\frac{p-q-1}{2}.
\end{equation}
Namely, (\ref{eq:pPsiqdefined}) converges absolutely for $|z|=\rho$ if $\Re(\mu)>0$, see \cite[Theorem~1]{KST}, \cite[Theorem~2.1]{KSST} and \cite[Theorem~1.5]{KSTBook}.  Note that $\mu$ in \eqref{eq:mu} is equal to $\mu-1/2$ in the notation of these references. The function ${}_{p}\Psi_q$ is a natural extension of the generalized hypergeometric function. It was introduced by Fox \cite{Fox} and Wright \cite{Wright35,Wright40,Wright52}, where both authors studied  its asymptotic behavior for large $z$ if $\Delta>-1$. The most well studied cases of ${_{p}\Psi_q}$ are the Wright (or the Bessel-Maitland) function ${_{0}\Psi_{1}}$  and the Mittag-Leffler function ${_{1}\Psi_{1}}$ with $A_1=a_1=1$.

The importance of the Fox-Wright function comes mostly from its role in fractional calculus, see \cite{GLM,Kilbas,KSTBook,MainardiPagnini}.  Other interesting applications also exist. Wright's original interest in this function was related to the asymptotic theory of partitions \cite[section~4.1]{Luchko}. In particular, Miller \cite{Miller} expressed a solution of the general trinomial equation in terms of ${}_1\!\Psi_{\!1}$. See also \cite[section~6]{MillerMosk} for an application in information theory.  It has been a recent surge of interest in the Fox-Wright function as witnessed by the articles \cite{ChuWang,Mehrez,MehrezSitnik,Paris2010,Paris2014,Paris2017,Prodanov,Wang,WeiGongHao}. The papers  \cite{ChuWang,WeiGongHao} establish summation formulas for the Fox-Wright function using combinatorial inversion formulas. In \cite{Wang} the author used a somewhat opposite approach by first developing contiguous relations for the Fox-Wright function and then employing them to  prove Hagen-Rothe convolutions from combinatorics. Paris \cite{Paris2010,Paris2014,Paris2017} discusses exponentially small asymptotics for ${_{p}\Psi_q}$ when $\Delta>-1$.  Mehrez \cite{Mehrez} extended our approach to generalized hypergeometric functions from \cite{KarpJMS2015,KLJAT2017,KPJMAA2012,KPITSF2017} to the Fox-Wright function ${}_{p}\Psi_q$ and presented  the Laplace and the generalized Stieltjes transform representations for some cases of ${}_{p}\Psi_q$ invoking certain facts from our papers \cite{KPCMFT2016,KPCMFT2017}.  He further obtained inequalities and monotonicity results  for  ${}_{p}\Psi_q$ established by us previously for the generalized hypergeometric function.  In another work \cite{MehrezSitnik} Mehrez and Sitnik derived further inequalities for the Fox-Wright functions and their ratios.
A recent survey of the Fox-Wright function and its applications can be found in \cite{Luchko}.

Before proceeding to the results of this paper it is convenient to introduce Fox's $H$ function.  Keeping earlier convention we assume that $\A=(A_1,\ldots,A_p)$ and  $\B=(B_1,\ldots,B_q)$ are positive scaling factors, $\a\in\C^p$ and $\b\in\C^{q}$, $0\le{n}\le{p}$, $1\le{m}\le{q}$.  Fox's $H$ function is defined by
\begin{equation}\label{eq:Fox}
H^{m,n}_{p,q}\left(z\left|\begin{array}{l} (\a,\A) \\(\b,\B) \end{array}\right.\right)
=\frac{1}{2\pi{i}}\int\limits_{\L}\frac{\prod\nolimits_{k=1}^{m}\Gamma(B_ks+b_k)\prod\nolimits_{i=1}^{n}\Gamma(1-a_i-A_is)}
{\prod\nolimits_{k=m+1}^{q}\Gamma(1-b_k-B_ks)\prod\nolimits_{i=n+1}^{p}\Gamma(A_is+a_i)}z^{-s}ds.
\end{equation}

The contour $\L$ separates the poles of the integrand at the points $s=-(b_k+j)/B_k$, $k\in\{1,\ldots,m\}$, $j\in\Z_{\ge0}$, from the poles at the points $s=(1-a_i+j)/A_i$, $i\in\{1,\ldots,n\}$, $j\in\Z_{\ge0}$.  It can be one of the following:
\begin{itemize}

\item $\L=\L_{-}$ is a left loop situated in a horizontal strip starting at the point
$-\infty+i\varphi_1$ and terminating at the point $-\infty+i\varphi_2$, where  $-\infty<\varphi_1<\varphi_2<\infty$;

\item $\L=\L_{+}$ is a right loop situated in a horizontal strip starting at the point
 $\infty+i\varphi_1$  and terminating at the point $\infty+i\varphi_2$, where $-\infty<\varphi_1<\varphi_2<\infty$;

\item $\L=\L_{i\gamma}$ coincides with the the vertical line $\Re{z}=\gamma$, $\gamma\in\R$, for sufficiently large $|z|$.

\end{itemize}
Details regarding the choice of the contour and the conditions for convergence of the integral in (\ref{eq:Fox}) can be found in \cite[sections 1.1,1.2]{KilSaig}, \cite[section~2.4]{ParKam} or \cite[section~8.3.1]{PBM3}.

In this paper we confine ourselves to the case $\Delta=-1$.  Under this restriction, the Fox-Wright function ${}_p\Psi_{q}$
 extends the  generalized hypergeometric function ${}_{p+1}F_{p}$ (the Gauss type) and the series definition  \eqref{eq:pPsiqdefined} only works for $|z|<\rho$  and $|z|=\rho$ if $\Re(\mu)>0$. Hence, we are encountered with the problem of analytic continuation of ${}_{p}\Psi_q$ beyond the disk of convergence of the series \eqref{eq:pPsiqdefined}.  This problem was completely solved by Braaksma \cite[Theorem~18]{Braaksma}, who derived his solution from more general results on Fox's $H$ function.

\newtheorem{thspecial}{Theorem}
\renewcommand\thethspecial{\Alph{thspecial}}

\begin{thspecial}\label{th:pPsiqanalytic}
For $\Delta=-1$ the function  ${}_{p}\Psi_q$ defined by the series \emph{\eqref{eq:pPsiqdefined}} can be continued analytically to $\C\setminus[\rho,\infty)$. The analytic continuation is given by

\emph{(1)} the Mellin-Barnes integral
\begin{equation}\label{eq:pPsiqMB}
{}_{p}\Psi_q(z)=\int_{\L_{i\gamma}}\frac{\Gamma(s)\Gamma(\a-\A{s})}{\Gamma(\b-\B{s})}(-z)^{-s}ds
\end{equation}
with some real $\gamma$ to the sector $|\arg(-z)|<\pi$\emph{;}

\emph{(2)} the convergent series
\begin{equation}\label{eq:pPsiqResidues}
{}_{p}\Psi_q(z)=\sum\limits_{\substack{{k=1,\ldots,p}\\{j=0,1,\ldots}}}\mathop{\mathrm{res}}\limits_{s=-(a_k+j)/A_k}\left\{\frac{\Gamma(s)\Gamma(\a-\A{s})}{\Gamma(\b-\B{s})}(-z)^{-s}\right\}
\end{equation}
from the sector $|\arg(-z)|<\pi$ to the domain $|z|>\rho$.
\end{thspecial}

If all poles of the integrand in \eqref{eq:pPsiqMB} are simple, then \eqref{eq:pPsiqResidues} takes the form \cite[(4.9)]{KSST}:
\begin{equation}\label{eq:pPsiqinv-series}
{}_{p}\Psi_q(z)=\sum\limits_{k=1}^{p}\sum\limits_{n=0}^{\infty}\frac{\Gamma(\a_{[k]}-\A_{[k]}(a_k+n)/A_k)}{\Gamma(\b-\B(a_k+n)/A_k)}
\Gamma\left(\frac{a_k+n}{A_k}\right)\frac{(-1)^n}{A_kn!}\left(-\frac{1}{z}\right)^{(a_k+n)/A_k},
\end{equation}
where the series converges for $|z|>\rho$ and $|z|=\rho$ if $\Re(\mu)>0$.

Note that representations and analytic continuation formulas for ${}_p\Psi_{q}(z)$ have been consequently considered in \cite{KST,KSST}, where the authors define the Fox-Wright function by integral \eqref{eq:pPsiqMB} over the contour $\L_{-}$.  The  authors also presented the series \eqref{eq:pPsiqinv-series} under this choice of contour, but did not make any claims regarding the analytic continuation of the series \eqref{eq:pPsiqdefined}.  It seems that Braaksma's results had been overlooked in these references.

Theorem~\ref{th:pPsiqanalytic} implies that the only singularity of ${}_{p}\Psi_q$ on the circle $|z|=\rho$ is the point $z=\rho$.  The behavior of ${}_{p}\Psi_q(z)$  in the neighborhood of this singular point and analytic continuation of ${}_{p}\Psi_q((\a,\A);(\b;\B);\rho)$ as the function of the parameters $\a,\A,\b,\B$ to a domain, where the series \eqref{eq:pPsiqdefined} diverges, seem to remain the only substantial unsolved problems in the theory of the Fox-Wright function (at least for $\Delta=-1$).  The main purpose of this paper is to present solutions of these two problems.  Their particular cases can be traced back to Gauss who found the celebrated summation formula for the hypergeometric function ${}_2F_{1}(a,b;c;1)$ furnishing analytic continuation in $a$, $b$, $c$, and described the behavior of $z\to{}_2F_{1}(a,b;c;z)$ in the neighborhood of $z=1$.  Similar problems for ${}_3F_{2}$ were partly solved by Thomae \cite{Thomae} in 1870 (analytic continuation in parameters) and partly by Ramanujan around the second decade of the 20th century (asymptotic approximation as $z\to1$ in the logarithmic case) with further contributions by Evans and Stanton, Wimp and B\"{u}hring, see \cite{Buehring87} and references therein. For the general Gauss type hypergeometric function ${}_{p+1}F_{p}$ with $p\ge3$ these problems were first solved by N{\o}rlund \cite{Norlund} with further contributions by  Olsson \cite{Olsson}, Marichev and Kalla \cite{MK}, Saigo and Srivastava \cite{SS}, B\"{u}hring \cite{Buehring92}, B\"{u}hring and  Srivastava \cite{BS} and other authors.

This paper is organized as follows.  In Section~2 we give a description of the behavior of ${}_{p}\Psi_q(z)$ in the neighborhood of $z=\rho$ in when the parameter $\mu$ from \eqref{eq:mu} is not an integer.  Section~3 treats the logarithmic cases of integer $\mu$. In Section~4 we deduce the formulas for the jump and the average value of ${}_{p}\Psi_q(z)$ when crossing the branch cut $[\rho,\infty)$.  Finally, Section~5 is devoted to conjectures and open problems.

\section{The Fox-Wright function near singularity}

We start with an integral representation of the Fox-Wright function which is a direct consequence of the Mellin transform formula for the delta-neutral $H$ function \cite[Theorem~6]{KPCMFT2016}.  First, define
\begin{equation}\label{eq:muprime}
\mu_{\sigma}=\mu+\sigma=\sum\nolimits_{j=1}^{q}b_j+\sigma-\sum\nolimits_{k=1}^{p}a_k+\frac{p-q-1}{2}.
\end{equation}
\begin{lemma}\label{lm:GSTFoxWright}
Suppose $\sigma>0$ is any number satisfying $\Re(\mu_{\sigma})>0$.
Assume further that $\Delta=-1$ and $\min_k(\Re(a_k/A_k))>0$.  Then for all $z\in\C\setminus[\rho,\infty)$
\begin{equation}\label{eq:GSTFoxWright}
{_{p}\Psi_q}\left(\left.\!\!\begin{array}{c}(\a,\A)\\(\b,\B)\end{array}\right|z\!\right)
=\Gamma(\sigma)\int_{0}^{1}\frac{H\!\left(\rho^{-1}{t}\right)dt}{t(1-\rho^{-1}{t}z)^{\sigma}},
\end{equation}
where
\begin{equation}\label{eq:HBprime}
H(z)=H_{q+1,p}^{p,0}\left(z\left|\begin{array}{l}(\b_{\sigma},\B_1)\\(\a,\A)\end{array}\right.\right)=\frac{1}{2\pi{i}}
\int\limits_{\L_{-}}\frac{\Gamma(\A{s}+\a)z^{-s}ds}{\Gamma(\B{s}+\b)\Gamma(s+\sigma)}
\end{equation}
and $\B_1=(B_1,\ldots,B_q,1)$, $\b_{\sigma}=(b_1,\ldots,b_q,\sigma)$.
\end{lemma}
\textbf{Proof.}  Note first that the equality $\Delta=-1$ is equivalent to $\sum_{j=1}^{q}B_j+1=\sum_{k=1}^{p}A_k$, so that $H(z)$ in \eqref{eq:GSTFoxWright} is delta-neutral in the terminology of \cite{KPCMFT2016,KPCMFT2017}. Then,  in view of the conditions $\Re(\mu_{\sigma})>0$ and $\min_k(\Re(a_k/A_k))>0$, we are in the position to apply \cite[Theorem~6]{KPCMFT2016} for the Mellin transform of $H_{q+1,p}^{p,0}$. Hence, we will have for $|z|<\rho$:
\begin{multline*}
\int_{0}^{1}\frac{H\!\left(\rho^{-1}{t}\right)dt}{t(1-\rho^{-1}{t}z)^{\sigma}}=\int_{0}^{\rho^{-1}}\frac{H(u)du}{u(1-uz)^{\sigma}}
=\sum\limits_{n=0}^{\infty}\frac{(\sigma)_n}{n!}z^n\int_{0}^{\rho^{-1}}u^{n-1}H(u)du
\\
=\sum\limits_{n=0}^{\infty}\frac{(\sigma)_n}{n!}z^n\frac{\Gamma(\A{n}+\a)}{\Gamma(\B{n}+\b)\Gamma(n+\sigma)}
=\frac{1}{\Gamma(\sigma)}{_{p}\Psi_q}\left(\left.\!\!\begin{array}{c}(\a,\A)\\(\b,\B)\end{array}\right|z\!\right).
\end{multline*}
To justify the term-wise integration we apply a version of the Lebesgue dominated convergence theorem  as presented in \cite[Theorem~3.6.2]{Bendetto}.  According to this theorem, it suffices to demonstrate that
$$
\sum\limits_{n=0}^{\infty}\frac{(\sigma)_n}{n!}|z|^n\int_{0}^{\rho^{-1}}u^{n-1}|H(u)|du
=\sum\limits_{n=0}^{\infty}\frac{(\sigma)_n}{n!}|z/\rho|^n\int_{0}^{1}t^{n-1}|H(\rho^{-1}t)|dt
<\infty.
$$
In view of the asymptotic behavior of $H\!\left(\rho^{-1}{t}\right)$ as $t\to0$ \cite[Corollary~1.12.1]{KilSaig} and $t\to1$ \cite[Theorem~1]{KPCMFT2017}, \cite[(24)]{Karp2019}, it follows from the assumptions $\Re(\mu_{\sigma})>0$ and $\min_k(\Re(a_k/A_k))>0$ that
$$
\int_{0}^{1}t^{n-1}|H\!\left(\rho^{-1}{t}\right)|dt<\int_{0}^{1}t^{-1}|H\!\left(\rho^{-1}{t}\right)|dt<M<\infty,
$$
since the integral on the left hand side is a decreasing function of $n$.  As
$$
\sum\limits_{n=0}^{\infty}\frac{(\sigma)_n}{n!}|z/\rho|^n<\infty
$$
for $|z|<\rho$ we are done.  The Fox-Wright function on the left hand side is analytic for $z\in\C\setminus[\rho,\infty)$ according to Theorem~\ref{th:pPsiqanalytic}, while the integral on the right is analytic in the same domain by uniform convergence on compact subsets.  Hence, by an analytic continuation, representation (\ref{eq:GSTFoxWright})  holds for all $z\in\C\setminus[\rho,\infty)$.   $\hfill\square$

Next, we will use the expansion for the delta-neutral Fox's $H$ function defined in \eqref{eq:HBprime} found in \cite{KPCMFT2017} and \cite{Karp2019}.  In order to present this expansion, we will need the Bernoulli-N{\o}rlund (or the generalized Bernoulli) polynomials $\Be^{(\sigma)}_{k}(x)$  defined by the generating function \cite[(1)]{Norlund61}:
$$
\frac{t^{\sigma}e^{xt}}{(e^t-1)^{\sigma}}=\sum\limits_{k=0}^{\infty}\Be^{(\sigma)}_{k}(x)\frac{t^k}{k!}.
$$
In particular, $\Be^{(1)}_{k}(x)=\Be_{k}(x)$ is the classical Bernoulli polynomial.  Further, the so-called signless non-central Stirling numbers of the first kind $s_{\sigma}(n,l)$ defined by their their ''horizontal'' generating function \cite[8.5]{Charalambides}
\begin{equation*}
(x+\sigma)_{n}=\sum\limits_{l=0}^{n}x^{l}s_{\sigma}(n,l),~~\text{where}~(a)_n=\Gamma(a+n)/\Gamma(a),
\end{equation*}
are related to the Bernoulli-N{\o}rlund polynomials by \cite[(7.6)]{Carlitz2}
$$
s_{\sigma}(n,l)=\binom{-l-1}{n-l}\Be^{(n+1)}_{n-l}(1-\sigma)=\frac{(-1)^{n-l}(l+1)_{n-l}}{(n-l)!}\Be^{(n+1)}_{n-l}(1-\sigma).
$$
Numerous formulas  for these numbers are collected in \cite{Charalambides}. Recall that $\B_1=(B_1,\ldots,B_q,1)$ and $\b_{\sigma}=(b_1,\ldots,b_q,\sigma)$. We can now define the coefficients
$l_r$, $r\ge1$, recursively, by
\begin{equation}\label{eq:lr}
l_r(\A,\a;\B_1,\b_{\sigma})=\frac{1}{r}\sum\limits_{m=1}^{r} q_m(\A,\a;\B_1,\b_{\sigma})l_{r-m}(\A,\a;\B_1,\b_{\sigma}),
\end{equation}
where $l_0=1$, and
\begin{equation}\label{eq:qm}
q_m(\A,\a;\B_1,\b_{\sigma})=\frac{(-1)^{m+1}}{m+1}\left[\sum\limits_{k=1}^p\frac{\Be_{m+1}(a_k)}{A_k^m}-\sum\limits_{j=1}^{q}\frac{\Be_{m+1}(b_j)}{B_j^m}-\Be_{m+1}(\sigma)\right].
\end{equation}
Similarly,  $l_0^{\theta}=1$, and
\begin{equation}\label{eq:lrtheta}
l_r^{\theta}(\A,\a;\B_1,\b_{\sigma})=\frac{1}{r}\sum\limits_{m=1}^r q_m^{\theta}(\A,\a;\B_1,\b_{\sigma}) l_{r-m}^{\theta}(\A,\a;\B_1,\b_{\sigma}),
\end{equation}
where
\begin{multline}\label{eq:qmtheta}
q_m^{\theta}(\A,\a;\B_1,\b_{\sigma})=\frac{(-1)^{m+1}}{m+1}
\\
\times\left[\sum\limits_{k=1}^p\frac{\Be_{m+1}(a_k)}{A_k^m}
-\sum\limits_{j=1}^q\frac{\Be_{m+1}(b_j)}{B_j^m}-\Be_{m+1}(\sigma)+\Be_{m+1}(\theta+\mu_{\sigma})-{\Be}_{m+1}(\theta+1)\right].
\end{multline}
Note that for $\mu_{\sigma}=1$ the numbers $l_r^{\theta}(\A,\a;\B_1,\b_{\sigma})$ do not depend on $\theta$ and reduce to $l_r(\A,\a;\B_1,\b_{\sigma})$.  We also remark that there are alternative ways to calculate $l_r$ from $q_m$ (and $l_r^{\theta}$ from $q_m^{\theta}$), namely
$$
l_{r}=\sum\limits_{\substack{{k_1+2k_2+\cdots+rk_r=r}\\{k_i\ge0}}}\frac{q_1^{k_1}(q_2/2)^{k_2}\cdots (q_r/r)^{k_r}}{k_1!k_2!\cdots k_r!}
=\sum\limits_{n=1}^{r}\frac{1}{n!}\sum\limits_{\substack{{k_1+k_2+\cdots+k_n=r}\\{k_i\ge1}}}\prod\limits_{i=1}^{n}\frac{q_{k_i}}{k_i}
$$
and
\begin{equation}\label{eq:lrBell}
l_{r}=\frac{1}{r!}Y_r(0!q_1,1!q_2,2!q_3,\ldots,(r-1)!q_r),
\end{equation}
where the (exponential) complete Bell polynomials are generated by \cite[(11.9)]{Charalambides}
$$
\exp\left(\sum\nolimits_{m=1}^{\infty}x_m\frac{t^m}{m!}\right)=1+\sum_{n=1}^{\infty}Y_n(x_1,\ldots,x_n)\frac{t^n}{n!},
$$
and given explicitly by \cite[(11.1)]{Charalambides}
$$
Y_n(x_1,\ldots,x_n)=\sum\limits_{\substack{{k_1+2k_2+\cdots+nk_n=n}\\{k_i\ge0}}}\frac{n!(x_1/1!)^{k_1}(x_2/2!)^{k_2}\cdots(x_n/n!)^{k_n}}{k_1!k_2!\cdots k_n!}.
$$
Furthermore, Nair found a determinantal expression for the solution of the recurrence (\ref{eq:lr}), which in our notation takes the form
$$
l_{r}=\frac{\det(\Omega_r)}{r!},~~~\Omega_r=[\omega_{i,j}]_{i,j=1}^{r},~~\omega_{i,j}\!=\!\left\{\!\!\begin{array}{ll}(i-1)!q_{i-j+1}/(j-1)!, & i\ge{j},\\
-1, &i=j-1,\\0, &i<j-1.
\end{array}\right.
$$
This formula could also be discovered by using the determinantal expression for the complete Bell polynomials \cite[p.203]{Collins}.  Define the constant $\nu$ by
\begin{equation}\label{eq:nu}
\nu=(2\pi)^{(p-q-1)/2}\prod\nolimits_{k=1}^{p}A_k^{a_k-1/2}\prod\nolimits_{j=1}^{q}B_j^{1/2-b_j}.
\end{equation}
The following theorem was proved in \cite[Theorem~1]{KPCMFT2017} with correction in \cite[Theorem~6]{Karp2019}.
\begin{theorem}\label{th:NorforH}
Suppose $\Delta=-1$, $\A,\B>1/6$, $\mu_{\sigma}=\mu+\sigma$ and $\theta$ is an arbitrary real number. Then $H_{q+1,p}^{p,0}\left(\rho^{-1}{t}\right)$ is represented by the series\emph{:}
\begin{equation}\label{eq:Hexpansion1}
H_{q+1,p}^{p,0}\left(\rho^{-1}{t}\left|\begin{array}{l}(\b_{\sigma},\B_1)\\(\a,\A)\end{array}\right.\right)
=t^{\theta+1}(1-t)^{\mu_{\sigma}-1}\sum\limits_{n=0}^{\infty}V_n(\theta)(1-t)^{n}
\end{equation}
convergent in the disk $|1-t|<1$, where $\B_1=(B_1,\ldots,B_q,1)$ and $\b_{\sigma}=(b_1,\ldots,b_q,\sigma)$,
\begin{equation}\label{eq:Vndefined}
V_n(\theta)\!=\!\frac{\nu}{\Gamma(n+\mu_{\sigma})}\sum\limits_{r=0}^{n}l_r^{\theta}(\A,\a;\B_1,\b_{\sigma})s_{\theta+\mu_{\sigma}}(n,r)
\!=\!\nu\!\sum\limits_{r+k=n}\!\frac{(-1)^kl_r(\A,\a;\B_1,\b_{\sigma})}{k!\Gamma(r+\mu_{\sigma})}\Be^{(n+\mu_{\sigma})}_{k}(-\theta)
\end{equation}
with $l_r(\A,\a;\B_1,\b_{\sigma})$ and $l_r^{\theta}(\A,\a;\B_1,\b_{\sigma})$ defined in \eqref{eq:lr} and \eqref{eq:lrtheta}, respectively.
\end{theorem}

We record the following estimate for the coefficients $V_n(\theta)$ to be used in the sequel.
\begin{lemma}\label{lm:Vnestimate}
Suppose $\Delta=-1$, $\A,\B>1/6$ and $\theta$ is arbitrary real number.  Let
$$
\alpha=\min_{1\le{j}\le{p}}\Re\left(a_j/A_j\right)
$$
and write $r$ for the maximal multiplicity among the poles of the integrand of $H_{q+1,p}^{p,0}\left(\rho^{-1}t\right)$ in  \eqref{eq:HBprime} with real part $-\alpha$.  Then for some constant $C>0$ and all $n\in\Z_{\ge2}$
\begin{equation}\label{eq:Vnestimate}
|V_n(\theta)|\le Cn^{\theta-\alpha}(\log{n})^{r-1}.
\end{equation}
\end{lemma}
\textbf{Proof.}  From \eqref{eq:Hexpansion1} we see that the function
$$
f_{\theta}(t)=t^{-\theta-1}(1-t)^{1-\mu_{\sigma}}H_{q+1,p}^{p,0}\left(\rho^{-1}{t}\left|\begin{array}{l}(\b_{\sigma},\B_1)\\(\a,\A)\end{array}\right.\right)
$$
is analytic in $|1-t|<1$ under conditions of the lemma.  Further, from \cite[Theorem~2]{Braaksma} (see details in \cite[Theorem~1]{Karp2019}) it follows that it is analytic in a larger domain
$$
\Delta(\phi,\varepsilon)=\{t:|1-t|\le 1+\varepsilon~\text~{and}~-\pi+\phi\le\arg(t)\le\pi-\phi\},
$$
for some $\varepsilon>0$ and any $0<\phi<\pi/2$.  According to \cite[(1.8.14)]{KilSaig} we have:
$$
f_{\theta}(t)=\mathcal{O}\left(t^{\alpha-\theta-1}(\log{t})^{r-1}\right)~\text{as}~t\to0.
$$
Then we are in the position to apply \cite[Theorem~2]{Flajolet} or \cite[Theorem~~V.3]{FSBook} to conclude that
$$
V_n=\mathcal{O}\left(n^{\theta-\alpha}(\log{n})^{r-1}\right),
$$
which implies \eqref{eq:Vnestimate}.$\hfill\square$

Our main result is the following theorem.
\begin{theorem}\label{th:FoxWrightrho}
Suppose $\Delta=-1$, $\A,\B>1/6$, $\alpha=\min_k(\Re(a_k/A_k))>0$ and $\mu\notin\Z$, where $\Delta$ and $\mu$ are defined in \eqref{eq:balance} and \eqref{eq:mu}, respectively. Then for $|1-z|<1/2$
\begin{equation}\label{eq:pPsiq2}
{_{p}\Psi_q}\left(\left.\!\!\begin{array}{c}(\a,\A)\\(\b,\B)\end{array}\right|\rho{z}\!\right)=
\left(1-z\right)^{\mu}\sum\limits_{m=0}^{\infty} R_m\left(1-z\right)^{m}
+\sum\limits_{m=0}^{\infty}W_m(1-z)^{m},
\end{equation}
where $\rho$ is given in \eqref{eq:rho}, and
\begin{align}
&R_m=-\frac{\pi\Gamma(m+\mu+\sigma)}{\sin(\pi\mu)}\sum\limits_{n=0}^{m}\frac{(-1)^{n}V_{n}(0)}
{\Gamma(\mu+n+1)(m-n)!}\label{eq:Rl},
\\
&W_m=(-1)^{m}\Gamma(\sigma+m)\sum\limits_{n=0}^{\infty}\frac{\Gamma(\mu+n-m)}{\Gamma(\mu+n+1)}V_n(0),
\label{eq:Wm}
\end{align}
and $V_n(0)$ are the coefficients given by \eqref{eq:Vndefined} with $\sigma$ being any positive number satisfying $\Re(\mu+\sigma)>0$.  The series in \eqref{eq:Wm} absolutely converges.
\end{theorem}

\noindent\textbf{Proof.}  Setting $\theta=0$ in  \eqref{eq:Hexpansion1} we obtain
\begin{equation}\label{eq:Hexpansion2}
H_{q+1,p}^{p,0}\left(\rho^{-1}{t}\left|\begin{array}{l}(\b_{\sigma},\B_1)\\(\a,\A)\end{array}\right.\right)
=t\sum\limits_{n=0}^{\infty}V_n(1-t)^{n+\mu_{\sigma}-1}.
\end{equation}
Substituting this expansion into  \eqref{eq:GSTFoxWright} and integrating termwise we will get
\begin{equation}\label{eq:pPsiq1}
\frac{1}{\Gamma(\sigma)}{_{p}\Psi_q}\left(\left.\!\!\begin{array}{c}(\a,\A)\\(\b,\B)\end{array}\right|\rho{z}\!\right)
=\sum\limits_{n=0}^{\infty}V_n\int_{0}^{1}\frac{(1-t)^{n+\mu_{\sigma}-1}dt}{(1-tz)^{\sigma}}
=\sum\limits_{n=0}^{\infty}\frac{V_n}{\mu_{\sigma}+n}
{}_{2}F_{1}\left.\!\!\left(\!\begin{matrix}1, \sigma\\\mu_{\sigma}+n+1\end{matrix}\right\vert z\!\right),
\end{equation}
where we applied the celebrated Euler's integral representation for the Gauss hypergeometric function ${}_2F_1$ \cite[Theorem~2.2.1]{AAR}. To justify the termwise integration we again apply the Lebesgue dominate convergence theorem in the form \cite[Theorem~3.6.2]{Bendetto}. According to this theorem the following estimate suffices
\begin{multline*}
\sum\limits_{n=0}^{\infty}|V_n|\int_{0}^{1}\frac{(1-t)^{n+\mu_{\sigma}-1}dt}{|1-zt|^{\sigma}}
\le\frac{1}{(M_z)^{\sigma}}\sum\limits_{n=0}^{\infty}\frac{|V_n|}{n+\mu_{\sigma}}
\\
\le\frac{C}{(M_z)^{\sigma}}\biggl\{\frac{|V_0|}{\mu_{\sigma}}+\frac{|V_1|}{\mu_{\sigma}+1}
+\sum\limits_{n=2}^{\infty}\frac{(\log{n})^{r-1}}{n^{\alpha}(n+\mu_{\sigma})}\biggr\}<\infty,
\end{multline*}
where we applied \eqref{eq:Vnestimate} and the last inequality holds as $\alpha>0$ by hypothesis.  Here $M_z$ is the distance between $1$ and the straight line segment $[0,z]$.

We can now use the connection formula \cite[Corollary~2.3.3]{AAR}
\begin{multline*}
{}_{2}F_{1}\left.\!\!\left(\begin{matrix}a, b\\c\end{matrix}\right\vert x\!\right)
\!
=\!\frac{\Gamma(c)\Gamma(a+b-c)}{\Gamma(a)\Gamma(b)}(1-x)^{c-a-b}{}_{2}F_{1}\left.\!\!\left(\!\begin{matrix}c-a, c-b\\c+1-a-b\end{matrix}\right\vert \!1-x\right)
\\
+\frac{\Gamma(c)\Gamma(c-a-b)}{\Gamma(c-b)\Gamma(c-a)}{}_{2}F_{1}\left.\!\!\left(\!\begin{matrix}a, b\\a+b+1-c\end{matrix}\right\vert 1-x\!\right)
\end{multline*}
valid if $c-a-b\notin\Z$ and $c\notin\Z_{\le0}$, to expand using $\mu_{\sigma}=\mu+\sigma$ and the binomial theorem:
\begin{multline}\label{eq:2F1series}
{}_{2}F_{1}\left.\!\!\left(\!\begin{matrix}1, \sigma\\\mu_{\sigma}+n+1\end{matrix}\right\vert z\!\right)
\!=\!\frac{\pi(-1)^{n+1}\Gamma(\mu_{\sigma}+n+1)(1-z)^{\mu+n}}{\Gamma(\sigma)\sin(\pi\mu)\Gamma(\mu+n+1)z^{\mu_{\sigma}+n}}
+\frac{\mu_{\sigma}+n}{\mu+n}{}_{2}F_{1}\left.\!\!\left(\!\begin{matrix}1, \sigma\\1-\mu-n\end{matrix}\right\vert1-z\!\right)
\\
=\frac{(1-z)^{\mu}(\mu_{\sigma}+n)}{\Gamma(\sigma)}\sum\limits_{m=0}^{\infty}G_{n,m}(1-z)^{m+n}
+\frac{\mu_{\sigma}+n}{\Gamma(\sigma)}\sum\limits_{m=0}^{\infty}D_{n,m}(1-z)^{m},
\end{multline}
where we applied the formula
$$
\Gamma(z+m)\Gamma(1-z-m)=(-1)^m\Gamma(z)\Gamma(1-z)=\frac{(-1)^m\pi}{\sin(\pi{z})}.
$$
The conditions $\mu\notin\Z$ and $\mu_{\sigma}+n+1\notin\Z_{\le0}$ must be satisfied for the validity of the above expansion. The latter is true due to the hypothesis $\Re(\mu+\sigma)>0$.  The coefficients take the form:
$$
G_{n,m}=\frac{(-1)^{n}\pi\Gamma(\mu_{\sigma}+n+m)}{\sin(\pi(\sigma-\mu_{\sigma}))\Gamma(\mu_{\sigma}-\sigma+n+1)m!}
=-\frac{(-1)^{n}\pi\Gamma(\mu_{\sigma}+n+m)}{\sin(\pi\mu)\Gamma(\mu+n+1)m!},
$$
$$
D_{n,m}=\frac{(-1)^{m}\Gamma(\sigma+m)\Gamma(\mu_{\sigma}+n-\sigma-m)}{\Gamma(\mu_{\sigma}+n-\sigma+1)}
=-\frac{\Gamma(\sigma+m)}{(-\mu-n)_{m+1}}.
$$
Substituting the first expansion in \eqref{eq:2F1series} into \eqref{eq:pPsiq1} we obtain:
\begin{multline}\label{eq:pPsi2F1}
\frac{1}{\Gamma(\sigma)}{_{p}\Psi_q}\left(\left.\!\!\begin{array}{c}(\a,\A)\\(\b,\B)\end{array}\right|\rho{z}\!\right)
=
-\frac{\pi(1-z)^{\mu}}{z^{\mu+\sigma}\sin(\pi\mu)\Gamma(\sigma)}
\sum\limits_{n=0}^{\infty}\frac{V_n\Gamma(\mu+n+\sigma)}{\Gamma(\mu+n+1)}\left(\frac{z-1}{z}\right)^{n}
\\
+\sum\limits_{n=0}^{\infty}\frac{V_n}{\mu+n}
{}_{2}F_{1}\left.\!\!\left(\!\begin{matrix}1, \sigma\\1-\mu-n\end{matrix}\right\vert 1-z\!\right).
\end{multline}

If, instead we substitute the second expansion  in \eqref{eq:2F1series} into \eqref{eq:pPsiq1} and interchange the order of summations, we will arrive at the representation:
\begin{equation*}
{_{p}\Psi_q}\left(\left.\!\!\begin{array}{c}(\a,\A)\\(\b,\B)\end{array}\right|\rho{z}\!\right)=
(1-z)^{\mu}\sum\limits_{m=0}^{\infty} \sum\limits_{n=0}^{\infty}V_{n}G_{n,m}(1-z)^{m+n}
+\sum\limits_{m=0}^{\infty}(1-z)^{m}\sum\limits_{n=0}^{\infty}V_{n}D_{n,m}.
\end{equation*}
Introducing the new index $n+m=l$ in the fist sum we can rewrite the above expansion as \eqref{eq:pPsiq2} with coefficients \eqref{eq:Rl} and \eqref{eq:Wm}.  To justify the interchange of the order of summations we will apply the Weierstrass double series theorem \cite[p.83]{Knopp}.  According to this theorem we need to prove that each  series in \eqref{eq:pPsi2F1} converges uniformly in $|1-z|\le{r}<1/2$. We will show that both series in fact converge uniformly in any compact subset of the domain $|(z-1)/z|<1$, which is easily seen to be the half plane $\Re(z)>1/2$; in particular, $|1-z|\le{r}<1/2$ is an example of such compact subset. Indeed, for the first series in \eqref{eq:pPsi2F1}  the claim follows immediately  from Lemma~\ref{lm:Vnestimate}.  The uniform convergence of the second series is seen from the asymptotic relation \cite[7.2(26)]{LukeBook}:
$$
{}_{2}F_{1}\left.\!\!\left(\!\begin{matrix}1, \sigma\\1-\mu-n\end{matrix}\right\vert 1-z\!\right)
\!=\!\left\{1+\frac{{\pi}n^{\sigma}z^{-\sigma}}{\Gamma(\sigma)\sin[\pi(n+\mu)]}
\left(\frac{1-z}{z}\right)^{n+\mu}\right\}\left(1+\mathcal{O}(n^{-1})\right),~n\to\infty,
$$
in view of Lemma~\ref{lm:Vnestimate} and the condition $\alpha>0$.  Then by the Weierstrass double series theorem, each series in  \eqref{eq:pPsiq2} converges for $|1-z|<1/2$ and the series in \eqref{eq:Wm} converges for each $m=0,1,\ldots$
The last claim is also implied by Lemma~\ref{lm:Vnestimate} and the asymptotic relation  (see, for instance, \cite[(1.2.1)]{KilSaig}):
$$
\frac{\Gamma(\mu+n-m)}{\Gamma(\mu+n+1)}=\frac{\Gamma(\mu_{\sigma}+n-\sigma-m)}{\Gamma(\mu_{\sigma}+n-\sigma+1)}\sim  \frac{1}{n^{m+1}},~~n\to\infty.~~~~~\square
$$

\textbf{Remark}. For arbitrary $\theta$ expansion \eqref{eq:pPsiq2} retains its form but the coefficients become
\begin{align*}
&R_m(\theta)=-\frac{\pi\Gamma(m+\mu+\sigma)\Gamma(m+\mu+\theta+1)}{\sin(\pi\mu)\Gamma(m+\mu+1)}
\sum\limits_{n=0}^{m}\frac{(-1)^{n}V_{n}(\theta)}{\Gamma(n+\mu+\theta+1)(m-n)!},
\\
&W_m(\theta)=\frac{\Gamma(m+\sigma)\Gamma(m+\theta+1)}{(-1)^{m}m!}
\sum\limits_{n=0}^{\infty}\frac{\Gamma(\mu+n-m)}{\Gamma(n+\mu+\theta+1)}V_n(\theta).
\end{align*}
Uniqueness of the power series coefficients imply that, in fact, $R_m(\theta)=R_m$, $W_m(\theta)=W_m$.  This independence of $\theta$ can be derived directly from the properties of the Bernoulli-N{\o}rlund polynomials derived by us in \cite[Theorem~6]{KPJMS2018} (\emph{check!}).

If $\Re(\mu)>-1$ we can choose $\sigma=1$ in \eqref{eq:GSTFoxWright}, so that formulas \eqref{eq:Rl} and \eqref{eq:Wm} simplify slightly. More interesting simplifications and transformations can be obtained by setting $\sigma=1$ in \eqref{eq:pPsi2F1}. This leads to the following
\begin{corollary}\label{cr:Foxfirstterm}
Suppose $\Re(\mu)>-1$ and all hypotheses of Theorem~\ref{th:FoxWrightrho} are satisfied. Then
\begin{align}
{_{p}\Psi_q}&\!\left(\left.\!\!\begin{array}{c}(\a,\A)\\(\b,\B)\end{array}\right|\rho{z}\!\right)
\!=\!-\frac{\pi(-1)^{\mu}}{\sin(\pi\mu)}H_{q+1,p}^{p,0}\left(\frac{1}{\rho{z}}\left|\begin{array}{l}(\b_{1},\B_1)\\(\a,\A)\end{array}\right.\!\!\right)
\!\!+\sum\limits_{n=0}^{\infty}\frac{V_n(0)}{\mu+n}
{}_{2}F_{1}\left.\!\!\left(\!\begin{matrix}1, 1\\1-\mu-n\end{matrix}\right\vert 1-z\!\right)\label{eq:pPsiq-SingH}
\\
&=-\frac{\pi(1-z)^{\mu}}{\sin(\pi\mu)z^{\mu}}\sum\limits_{m=0}^{\infty}(1-z)^m\Delta^{m}V_0(0)
+\sum\limits_{m=0}^{\infty}W_m(1-z)^{m},\label{eq:pPsiq-SingD}
\end{align}
where $V_n(0)$ and $W_m$ are the coefficients given in \eqref{eq:Vndefined}, \eqref{eq:Wm}, respectively, with $\sigma=1$ and
$$
\Delta^{m}V_0(0)=\sum\limits_{j=0}^{m}\binom{m}{j}(-1)^jV_j(0)
$$
is $m$-th forward difference of the sequence $\{V_j(0)\}_{j=1}^{\infty}$.
\end{corollary}

\textbf{Remark.} The Fox's $H$ function $H_{q+1,p}^{p,0}$ in \eqref{eq:pPsiq-SingH} for $|z|>1$ is understood as analytic continuation of the function defined in \eqref{eq:Fox}. This analytic continuation is given by the series \eqref{eq:Hexpansion1} or by contour integral as explained in \cite[Theorem~2]{Karp2019}.

\textbf{Proof.} Indeed, on setting $\sigma=1$ formula \eqref{eq:pPsi2F1} takes the form:
$$
{_{p}\Psi_q}\left(\left.\!\!\begin{array}{c}(\a,\A)\\(\b,\B)\end{array}\right|\rho{z}\!\right)
=-\frac{\pi(1-z)^{\mu}}{z^{\mu+1}\sin(\pi\mu)}\sum\limits_{n=0}^{\infty}V_n\left(\frac{z-1}{z}\right)^{n}
+\sum\limits_{n=0}^{\infty}\frac{V_n}{\mu+n}
{}_{2}F_{1}\left.\!\!\left(\!\begin{matrix}1, 1\\1-\mu-n\end{matrix}\right\vert 1-z\!\right).
$$
Comparing this formula with \eqref{eq:Hexpansion1} we immediately get \eqref{eq:pPsiq-SingH}.
The second term in \eqref{eq:pPsiq-SingD} is the same as in \eqref{eq:pPsiq2}. The first term in \eqref{eq:pPsiq-SingD}
is a result of an application of Euler's transformation \cite[(1.20)]{Norlund61}
$$
\sum\limits_{n=0}^{\infty}\frac{(\alpha)_n}{n!}f_nz^{n}
=(1-z)^{-\alpha}\sum\limits_{n=0}^{\infty}\frac{(\alpha)_n}{n!}\Delta^{n}f_0\left(\frac{z}{1-z}\right)^n
$$
with $\alpha=1$.  Here $\Delta^{n}f_0$ is the $n$-th difference, $\Delta{f_0}=f_1-f_0$, $\Delta^n{f_0}=\Delta(\Delta^{n-1}{f_0})$.$\hfill\square$

It seems worthwhile to write the principal term in \eqref{eq:pPsiq2} in a more explicit form.  In view of
$\Be^{\alpha}_{0}(x)=1$, formulas \eqref{eq:Rl}, \eqref{eq:Wm} and \eqref{eq:Vndefined} lead immediately to

\begin{corollary}\label{cr:pPsiq-asymp}
Under conditions of Theorem~\ref{th:FoxWrightrho} the following asymptotic relation is true as $z\to1$\emph{:}
$$
{_{p}\Psi_q}\left(\left.\!\!\begin{array}{c}(\a,\A)\\(\b,\B)\end{array}\right|\rho{z}\!\right)=
\nu\Gamma(-\mu)\left(1-z\right)^{\mu}(1+\mathcal{O}(1-z))+\Gamma(\sigma)\big(1+\mathcal{O}(1-z))\big)\sum\limits_{n=0}^{\infty}\frac{V_n(0)}{\mu+n},
$$
where $\mu$ and $\nu$ is defined by \eqref{eq:mu} and \eqref{eq:nu}, respectively; $V_n(0)$ are the coefficients given in \eqref{eq:Vndefined}.  In particular, if $\Re(\mu)>0$, we have the summation formula
\begin{equation}\label{eq:Psi-rho}
{_{p}\Psi_q}\left(\left.\!\!\begin{array}{c}(\a,\A)\\(\b,\B)\end{array}\right|\rho\!\right)=\Gamma(\sigma)\sum\limits_{n=0}^{\infty}\frac{V_n(0)}{\mu+n}.
\end{equation}
If $\sigma$ is chosen to satisfy $\sigma>\max(-\Re(\mu),0)$, the series on the right hand side converges for all $\A,\B>1/6$ and $\min_k(\Re(a_k/A_k))>0$ regardless of the value of $\mu$. Thus it provides the analytic continuation in parameters for the function on the left hand side.
\end{corollary}
\textbf{Remark.} Formula \eqref{eq:Psi-rho} can be viewed as a far-reaching generalization of the Gauss summation formulas for ${}_2F_{1}$ and further summation formulas for the generalized hypergeometric function ${}_{p+1}F_{p}(1)$ considered in \cite{Buehring92,Norlund,Olsson,SS}.

\section{The logarithmic cases}
Theorem~\ref{th:FoxWrightrho} fails if $\mu\in\Z$. The underlying reason for this is the failure of representation \eqref{eq:2F1series}. According to \cite[2.10(12)]{HTF1} for $\mu+n=s\in\Z_{\ge0}$ and $|\arg(1-z)|<\pi$ it takes the form:
\begin{multline}\label{eq:2F1logpositive}
\frac{1}{\Gamma(\mu_{\sigma}+n+1)}{}_{2}F_{1}\left.\!\!\left(\!\begin{matrix}1, \sigma\\\mu_{\sigma}+n+1\end{matrix}\right\vert z\!\right)
\!=\!\frac{1}{s\Gamma(\sigma+s)}\sum_{m=0}^{s-1}\frac{(\sigma)_m(1-z)^m}{(1-s)_m}
\\
+\frac{(-1)^s(1-z)^s}{s!\Gamma(\sigma)}\sum_{m=0}^{\infty}\frac{(\sigma+s)_m}{m!}(1-z)^m\big[-\log(1-z)+h_m\big]
\\
=-\sum_{m=0}^{s-1}\frac{(\sigma)_m(1-z)^m}{\Gamma(\sigma+s)(-s)_{m+1}}
-\frac{(-1)^s(1-z)^s}{z^{s+\sigma}s!\Gamma(\sigma)}\log(1-z)
+
\frac{(-1)^s(1-z)^s}{s!\Gamma(\sigma)}\sum_{m=0}^{\infty}h_m\frac{(\sigma+s)_m}{m!}(1-z)^m,
\end{multline}
where
$$
h_{s,m}=\psi(m+1)-\psi(m+s+\sigma),
$$
and $\psi(z)=\Gamma'(z)/\Gamma(z)$ denotes the digamma function.  For $\mu+n=s\in\Z_{<0}$  by  \cite[2.10(14)]{HTF1}:
\begin{equation}\label{eq:2F1lognegative}
\frac{1}{\Gamma(\mu_{\sigma}+n+1)}
{}_{2}F_{1}\left.\!\!\left(\!\begin{matrix}1,\sigma\\\mu_{\sigma}+n+1\end{matrix}\right\vert z\!\right)
\!=\!\frac{\Gamma(-s)(1-z)^{s}}{\Gamma(\sigma)}\sum_{m=0}^{-s-1}\frac{(\sigma+s)_m}{m!}(1-z)^m.
\end{equation}
Using these expansions we arrive at the following theorem.

\begin{theorem}\label{th:FoxWrightlogcase}
Suppose $\Delta=-1$, $\A,\B>1/6$, $\alpha=\min_k(\Re(a_k/A_k))>0$ and $\mu\in\Z_{\ge0}$, where $\Delta$ and $\mu$ are defined in \eqref{eq:balance} and \eqref{eq:mu}, respectively. Then for $|1-z|<1/2$
\begin{multline}\label{eq:pPsiqlogpos}
{_{p}\Psi_q}\left(\left.\!\!\begin{array}{c}(\a,\A)\\(\b,\B)\end{array}\right|\rho{z}\!\right)
=(1-z)^{\mu}\sum\limits_{j=0}^{\infty}\Gamma(\mu_{\sigma}+j)(1-z)^{j}
\sum\limits_{n+m=j}\frac{(-1)^{\mu+n}V_n}{(\mu+n)!m!}\bigl\{h_{n+\mu,m}-\log(1-z)\bigr\}
\\
-(1-z)^{\mu}\sum\limits_{j=-\mu}^{\infty}\Gamma(\mu_{\sigma}+j)(1-z)^{j}
\sum_{n=\max(0,j+1)}^{\infty}\frac{V_n}{(-\mu-n)_{\mu+j+1}},
\end{multline}
where $\rho$ is given in \eqref{eq:rho}, $V_n=V_n(0)$ are the coefficients given by \eqref{eq:Vndefined} with any $\sigma>0$ and $h_{n+\mu,m}=\psi(m+1)-\psi(m+n+\mu+\sigma)$.

If $\mu\in\Z_{<0}$, $\sigma>-\mu$, we have for $|1-z|<1/2$\emph{:}
\begin{multline}\label{eq:pPsiqlognegative}
{_{p}\Psi_q}\left(\left.\!\!\begin{array}{c}(\a,\A)\\(\b,\B)\end{array}\right|\rho{z}\!\right)
=\sum\limits_{j=0}^{\infty}\Gamma(\sigma+j)(1-z)^{j}
\sum\limits_{n+m=j}\frac{(-1)^{n}V_{n-\mu}}{n!m!}\bigl\{h_{n,m}-\log(1-z)\bigr\}
\\
-\sum\limits_{j=0}^{\infty}\Gamma(\sigma+j)(1-z)^{j}\sum_{n=j+1}^{\infty}\frac{V_{n-\mu}}{(-n)_{j+1}}
+\sum\limits_{j=\mu}^{-1}\Gamma(\sigma+j)(1-z)^j\sum\limits_{n=0}^{j-\mu}V_n\frac{(-\mu-n-1)!}{(-\mu-n+j)!}.
\end{multline}
\end{theorem}

\textbf{Proof.} Suppose first that $\mu\in\Z_{\ge0}$. Then also $\mu+n=s\in\Z_{\ge0}$ and substitution of \eqref{eq:2F1logpositive} into \eqref{eq:pPsiq1} formally yields ($V_n:=V_n(0)$ for brevity):
\begin{multline}\label{eq:pPsiqlog0}
{_{p}\Psi_q}\left(\left.\!\!\begin{array}{c}(\a,\A)\\(\b,\B)\end{array}\right|\rho{z}\!\right)
=\sum\limits_{n=0}^{\infty}V_n\Gamma(\mu_{\sigma}+n)\biggl\{-\sum_{m=0}^{\mu+n-1}\frac{\Gamma(\sigma+m)(1-z)^m}{\Gamma(\mu_{\sigma}+n)(-\mu-n)_{m+1}}
\\
+\frac{(-1)^{\mu+n}}{(\mu+n)!}\sum_{m=0}^{\infty}h_{n+\mu,m}\frac{(\mu_{\sigma}+n)_m}{m!}(1-z)^{\mu+n+m}
-\frac{(-1)^{\mu+n}(1-z)^{\mu+n}}{z^{\mu+n+\sigma}(\mu+n)!}\log(1-z)\biggr\}.
\end{multline}
Expanding $z^{-(\mu+n+\sigma)}$ in the last term in powers of $1-z$ by the binomial theorem and interchanging the order of summations after some simplifications we get \eqref{eq:pPsiqlogpos}.

To justify the interchange of the order of summations we estimate each term in  \eqref{eq:pPsiqlog0} separately.  First note that for $\mu+n=s\in\Z_{\ge0}$ and $m=0,\ldots,s-1$ we have
$$
|(-s)_{m+1}|\ge sm!
$$
so that, in view of Lemma~\ref{lm:Vnestimate}, for $|1-z|<1$:
\begin{multline*}
\sum\limits_{s=\mu}^{\infty}|V_{s-\mu}|\sum_{m=0}^{s-1}\frac{\Gamma(\sigma+m)|1-z|^m}{|(-s)_{m+1}|}
\le\Gamma(\sigma)\sum\limits_{s=\mu}^{\infty}\frac{|V_{s-\mu}|}{s}\sum_{m=0}^{\infty}\frac{(\sigma)_m|1-z|^m}{m!}
\\
\le\frac{\Gamma(\sigma)}{(1-|1-z|)^{\sigma}}\sum\limits_{s=\mu}^{\infty}\frac{|V_{s-\mu}|}{s}<\infty.
\end{multline*}
This justifies the interchange of the order of summations for the first term in braces in \eqref{eq:pPsiqlog0} by Tonelli's theorem.  To handle the second terms we can employ the bound \cite[(2.2)]{Alzer}
$$
\log(x)-\frac{1}{x}<\psi(x)<\log(x)-\frac{1}{2x},~~~x>0,
$$
which implies that ($s=n+\mu$)
$$
|h_{s,m}|=-h_{s,m}=\psi(m+s+\sigma)-\psi(m+1)\le\log\left(\frac{m+s+\sigma}{m+1}\right)+\frac{1}{m+1}-\frac{1}{m+s+\sigma}.
$$
According to the classical inequality $\log(1+x)\le{x}$ this leads to
$$
|h_{s,m}|\le\frac{s+\sigma}{m+1}-\frac{1}{m+s+\sigma}<\frac{s+\sigma}{m+1}\le \frac{s+\sigma+m}{m+1}.
$$
Hence, writing $x=|1-z|$ for the absolute value of the second term we get the bound
\begin{multline*}
\sum\limits_{s=\mu}^{\infty}|V_{s-\mu}|\Gamma(s+\sigma)
\frac{x^s}{s!}\sum_{m=0}^{\infty}|h_{s,m}|\frac{(s+\sigma)_m}{m!}x^{m}
\leq
\sum\limits_{s=\mu}^{\infty}|V_{s-\mu}|\Gamma(s+\sigma)
\frac{x^{s-1}}{s!}\sum_{m=0}^{\infty}\frac{(s+\sigma)_{m+1}}{(m+1)!}x^{m+1}
\\
\leq\frac{1}{x(1-x)^{\sigma}}\sum\limits_{s=\mu}^{\infty}|V_{s-\mu}|\frac{\Gamma(s+\sigma)}{\Gamma(s+1)}(x/(1-x))^s<\infty.
\end{multline*}
As for $x=|1-z|<1/2$ we have $x/(1-x)<1$ the last inequality follows by Lemma~\ref{lm:Vnestimate}.  Thus, Tonelli's theorem justifies the interchange of the order of summations.  The term containing logarithm is precisely the same as the term just considered but simpler ($h_{s,m}=1$) so that the above argument also covers this term.

For $\mu\in\Z_{<0}$ we split the summation in \eqref{eq:pPsiq1} in two parts: $n=0,\ldots-\mu-1$ and $n\ge-\mu$.
Then we apply \eqref{eq:2F1lognegative} for ${}_2F_1$ in the first part, and \eqref{eq:2F1logpositive} in the second.
As we only have a finite number of terms in the first part while the second part is similar to the case $\mu\in\Z_{\ge0}$, the justification of the interchange of the order of summations is not different from the one given above.  Performing this interchange after some simplifications we arrive at  \eqref{eq:pPsiqlognegative}.$\hfill\square$

In particular, if $\mu=0$ by choosing $\sigma=1$ (note that the condition $\mu+\sigma>0$ is then satisfied) we will have an expansion for what can be designated as the ''zero-balanced case'' of the Fox-Wright function.

\begin{corollary}\label{cr:pPsiq-zerobal}
Suppose $\mu=0$  and conditions of Theorem~\ref{th:FoxWrightlogcase} are satisfied. Then
$$
{_{p}\Psi_q}\left(\left.\!\!\begin{array}{c}(\a,\A)\\(\b,\B)\end{array}\right|\rho{z}\!\right)
=\sum\limits_{j=0}^{\infty}j!(1-z)^{j}\left[
\sum\limits_{n=0}^{j}\frac{(-1)^{n}V_n}{n!(j-n)!}\bigl\{h_{n,j-n}-\log(1-z)\bigr\}
-\sum_{n=j+1}^{\infty}\frac{V_n}{(-n)_{j+1}}\right],
$$
where $\rho$ is given in \eqref{eq:rho}, $V_n=V_n(0)$ are the coefficients given by \eqref{eq:Vndefined} with $\sigma=1$ and $h_{n,j-n}=\psi(j-n+1)-\psi(j+1)$.
In particular, the following asymptotic formula is true as $z\to1$\emph{:}
$$
{_{p}\Psi_q}\left(\left.\!\!\begin{array}{c}(\a,\A)\\(\b,\B)\end{array}\right|\rho{z}\!\right)
=-V_0\log(1-z)+\sum_{n=1}^{\infty}\frac{V_n}{n}+\mathcal{O}((1-z)\log(1-z)).
$$
\end{corollary}

\section{The behavior on the branch cut}
In this section, we will compute the jump and the average value of the Fox-Wright function on the branch cut $[\rho,\infty)$.
\begin{theorem}\label{th:Psi-jump}
Suppose $x>1$ and  $\a,\b$ are real vectors.  Then the following identities hold true
\begin{equation}\label{eq:Psi-jump}
{_{p}\Psi_q}\!\left(\!\!\begin{array}{l}(\a,\A)\\(\b,\B)\end{array}\vline\,\,x+i0\!\right)
-{_{p}\Psi_q}\!\left(\!\!\begin{array}{l}(\a,\A)\\(\b,\B)\end{array}\vline\,\,x-i0\!\right)
\!=\!2\pi{i}H_{q+1,p}^{p,0}\left(\frac{1}{x}\left|\begin{array}{l}(\b_{1},\B_1)\\(\a,\A)\end{array}\right.\!\!\right)
\end{equation}
and
\begin{equation}\label{eq:Psi-aver}
{_{p}\Psi_q}\!\left(\!\!\begin{array}{l}(\a,\A)\\(\b,\B)\end{array}\vline\,\,x+i0\!\right)
+{_{p}\Psi_q}\!\left(\!\!\begin{array}{l}(\a,\A)\\(\b,\B)\end{array}\vline\,\,x-i0\!\right)
\!=\!-2{\pi}H_{q+2,p+1}^{p,1}\left(\frac{1}{x}\left|\begin{array}{l}(1,1),(3/2,1),(\b,\B)\\(\a,\A),(3/2,1)\end{array}\right.\!\!\right).
\end{equation}
\end{theorem}
\textbf{Proof.}  If $\mu>-1$  and $H\!\left(\rho^{-1}{t}\right)\ge0$ for $t\in(0,1)$ we can choose $\sigma=1$ in \eqref{eq:GSTFoxWright}, so that
$$
{_{p}\Psi_q}\left(\left.\!\!\begin{array}{c}(\a,\A)\\(\b,\B)\end{array}\right|z\!\right)
=\int_{0}^{1}\frac{H\!\left(\rho^{-1}{t}\right)dt}{t(1-\rho^{-1}{t}z)}=\int_{\rho}^{\infty}\frac{\lambda(x)dx}{x+z},
$$
where
$$
\lambda(x)=H_{q+1,p}^{p,0}\!\left(\frac{1}{x}\left|\begin{array}{l}(\b_{1},\B_1)\\(\a,\A)\end{array}\right.\!\!\right).
$$
According to the Stieltjes inversion formula \cite[Theorem~7a,p.339]{Widder}
$$
\lambda(x)=\frac{1}{2\pi{i}}\left\{{_{p}\Psi_q}\left(\left.\!\!\begin{array}{c}(\a,\A)\\(\b,\B)\end{array}\right|x+i0\right)
-{_{p}\Psi_q}\left(\left.\!\!\begin{array}{c}(\a,\A)\\(\b,\B)\end{array}\right|x-i0\right)\right\},
$$
which gives formula \eqref{eq:Psi-jump} under the above restrictions. To remove these restrictions we assume that
all poles of the integrand in \eqref{eq:pPsiqMB} are simple so that ${_{p}\Psi_q}(z)$ can be represented by the convergent series \eqref{eq:pPsiqinv-series}.  Writing
$$
D_{n,k}=\frac{\Gamma(\a_{[k]}-\A_{[k]}(a_k+n)/A_k)}{\Gamma(\b-\B(a_k+n)/A_k)}\frac{(-1)^n}{A_kn!}
$$
and using $(-x-i0)^{-u}-(-x+i0)^{-u}={2i}x^{-u}\sin(\pi{u})$ for $x>1$ and real $u$, we calculate
\begin{multline*}
{_{p}\Psi_q}(x+i0)-{_{p}\Psi_q}(x-i0)\!
\\
=\!\sum\limits_{k=1}^{p}\sum\limits_{n=0}^{\infty}
D_{n,k}\Gamma\!\left(\frac{a_k+n}{A_k}\right)\left[(-x-i0)^{-(a_k+n)/A_k}-(-x+i0)^{-(a_k+n)/A_k}\right]
\\
=\!2i\!\sum\limits_{k=1}^{p}\sum\limits_{n=0}^{\infty}
D_{n,k}\Gamma\!\left(\frac{a_k+n}{A_k}\right)(1/x)^{(a_k+n)/A_k}\sin\!\left(\pi\frac{a_k+n}{A_k}\right)
\!
\\
=\!2\pi{i}\sum\limits_{k=1}^{p}\sum\limits_{n=0}^{\infty}\frac{D_{n,k}(1/x)^{(a_k+n)/A_k}}{\Gamma\left(1-(a_k+n)/A_k\right)},
\end{multline*}
where we omitted the parameter list in the notation of ${_{p}\Psi_q}$ for brevity (it is seen in \eqref{eq:Psi-jump}).
To see that the last series equals the right hand side of \eqref{eq:Psi-jump} it remains to apply \cite[Theorem~1.3]{KilSaig}. Finally, if the poles  of the integrand in \eqref{eq:pPsiqMB} are not simple \eqref{eq:Psi-jump} remains true by continuity of both sides.

For the average value we have in a similar fashion, but this time using
$$
(-x-i0)^{-u}+(-x+i0)^{-u}={2}x^{-u}\cos\pi{u}=-\frac{{2}x^{-u}\pi}{\Gamma(u-1/2)\Gamma(3/2-u)},
$$
the following chain
\begin{multline*}
{_{p}\Psi_q}(x+i0)+{_{p}\Psi_q}(x-i0)
\\
=\sum\limits_{k=1}^{p}\sum\limits_{n=0}^{\infty}
D_{n,k}\Gamma\left(\frac{a_k+n}{A_k}\right)\left\{(-x-i0)^{-(a_k+n)/A_k}+(-x+i0)^{-(a_k+n)/A_k}\right\}
\\
=-2\pi\sum\limits_{k=1}^{p}\sum\limits_{n=0}^{\infty}
\frac{D_{n,k}\Gamma\left((a_k+n)/A_k\right)(1/x)^{(a_k+n)/A_k}}{\Gamma(3/2-(a_k+n)/A_k)\Gamma(-1/2+(a_k+n)/A_k)}.
\end{multline*}
Another application of \cite[Theorem~1.3]{KilSaig} shows that this series is equal to the right hand side of
\eqref{eq:Psi-aver}.$\hfill\square$

\section{Conjectures and open problems}

Convergence of the series \eqref{eq:pPsiq2} in Theorem~\ref{th:FoxWrightrho} and the series \eqref{eq:pPsiqlogpos}, \eqref{eq:pPsiqlognegative} in Theorem~\ref{th:FoxWrightlogcase} in the disk $|1-z|<1/2$ results from the estimate of the radius of convergence of these series \emph{before} in the interchange of the order of summations.  The number $1/2$ is not claimed to be the true radius of convergence of these series. If it were so for one term it would necessarily mean that both terms in, say  \eqref{eq:pPsiq2}, have singularities on the circle  $|1-z|=1/2$ that cancel out, since the Fox-Wright function on the left hand side does not have such singular points. While this is possible in principle, such situation seems rather unlikely to be true.  This motivates the following conjecture.

\begin{conjecture} The true radius of convergence of all series in \eqref{eq:pPsiq2}, \eqref{eq:pPsiqlogpos} and \eqref{eq:pPsiqlognegative} is equal to unity.
\end{conjecture}

An important case not covered in this paper is $0<\min(\A,\B)<1/6$.  In this case $H^{p,0}_{q+1,p}$ in Theorem~\ref{th:NorforH} has singularities in the disk $|1-t|<1$ at the distance  $2\sin(\gamma_1/2)\in(0,1)$, $\gamma_1=2\pi\min(\A,\B)$, from the point $t=1$, so that the radius of convergence in \eqref{eq:Vndefined} is less than $1$. Details can be found in \cite{Karp2019}. This leads to a faster growth of the coefficients $V_n(\theta)$ in \eqref{eq:Vndefined} and divergence in \eqref{eq:Wm}.
To tackle this problem we can apply the relation \cite[(2.1.4)]{KilSaig}
\begin{equation*}\label{eq:Hrescaling}
H^{p,0}_{q+1,p}\left(\!\rho^{-1}{u^{1/\omega}}\left|\begin{array}{l}(\b_{\sigma},\B_1)\\(\a,\A)\end{array}\right.\!\right)
={\omega}H^{p,0}_{q+1,p}\left(\!\rho^{-\omega}u\left|\begin{array}{l}(\b_{\sigma},\omega\B_1)\\(\a,\omega\A)\end{array}\right.\!\right)
\end{equation*}
valid for any $\omega>0$.  Making substitution $t=u^{1/\omega}$ in \eqref{eq:GSTFoxWright} and applying the above relation, we obtain
\begin{equation*}\label{eq:Psi-int-omega}
\frac{1}{\Gamma(\sigma)}{_{p}\Psi_q}\left(\left.\!\!\begin{array}{c}(\a,\A)\\(\b,\B)\end{array}\right|z\!\right)
=\frac{1}{\omega}\int_{0}^{1}\frac{H\!\left(\rho^{-1}u^{1/\omega}\right)du}{u(1-\rho^{-1}u^{1/\omega}z)^{\sigma}}
=\int_{0}^{1}\frac{\hat{H}\!\left(\rho^{-\omega}u\right)du}{u(1-\rho^{-1}u^{1/\omega}z)^{\sigma}},
\end{equation*}
where by Theorem~\ref{th:NorforH} (note that by \eqref{eq:rho} the number $\rho^{-\omega}$ plays the role of $\rho^{-1}$ under substitution $\A\to\omega\A$, $\B_1\to\omega\B_1$)
\begin{equation*}\label{eq:hatH}
\hat{H}\!\left(\rho^{-\omega}u\right)=H^{p,0}_{q+1,p}\left(\rho^{-\omega}u\left|\begin{array}{l}(b_{\sigma},\omega\B_1)\\(\a,\omega\A)\end{array}\right.\!\right)
=u^{\theta+1}(1-u)^{\mu_{\sigma}-1}\sum\limits_{n=0}^{\infty}\hat{V}_n(\theta)(1-u)^{n}.
\end{equation*}
Here the coefficients $\hat{V}_n(\theta)$ are still computed by \eqref{eq:Vndefined} but taking $l^{\theta}_r(\omega\A,\a;\omega\B_1,\b_{\sigma})$ and $l_r(\omega\A,\a;\omega\B_1,\b_{\sigma})$ instead of
$l^{\theta}_r(\A,\a;\B_1,\b_{\sigma})$ and $l_r(\A,\a;\B_1,\b_{\sigma})$, respectively. Choosing $\omega>1$ large enough to satisfy the condition $\omega\A,\omega\B>1/6$ we guarantee the convergence of the above expansion in the disk $|1-u|<1$ and validity of the argument given in the proof of Theorem~\ref{th:FoxWrightrho} to justify the term-wise integration.  Hence, taking as before $\theta=0$ we get
$$
\frac{1}{\Gamma(\sigma)}{_{p}\Psi_q}\!\left(\left.\!\!\begin{array}{c}(\a,\A)\\(\b,\B)\end{array}\right|\rho{z}\!\right)
=\sum\limits_{n=0}^{\infty}\hat{V}_n(0)\int_{0}^{1}\frac{(1-u)^{n+\mu_{\sigma}-1}du}{(1-u^{1/\omega}z)^{\sigma}}.
$$
The integral on the right can be computed by the standard termwise integration:
\begin{multline}\label{eq:2Psi1int}
\int_{0}^{1}\frac{(1-u)^{n+\mu_{\sigma}-1}du}{(1-u^{1/\omega}z)^{\sigma}}
=\sum\limits_{j=0}^{\infty}\frac{(\sigma)_j}{j!}z^{j}\int_{0}^{1}u^{j/\omega}(1-u)^{n+\mu_{\sigma}-1}du
\\
=\sum\limits_{j=0}^{\infty}\frac{(\sigma)_j}{j!}z^{j}\frac{\Gamma(j/\omega+1)}{\Gamma(j/\omega+n+\mu_{\sigma}+1)}
=\frac{\Gamma(n+\mu_{\sigma})}{\Gamma(\sigma)}
{_{2}\Psi_1}\!\left(\left.\!\!\begin{array}{c}(\sigma,1),(1,1/\omega)\\(n+\mu_{\sigma}+1,1/\omega)\end{array}\right|z\!\right).
\end{multline}
Hence,  in place of formula \eqref{eq:pPsiq1} we now obtain:
\begin{equation}\label{eq:pPsiq2Psi1}
{_{p}\Psi_q}\!\left(\left.\!\!\begin{array}{c}(\a,\A)\\(\b,\B)\end{array}\right|\rho{z}\!\right)
=\sum\limits_{n=0}^{\infty}\hat{V}_n(0)
\Gamma(n+\mu_{\sigma})
{_{2}\Psi_1}\!\left(\left.\!\!\begin{array}{c}(\sigma,1),(1,1/\omega)\\(n+\mu_{\sigma}+1,1/\omega)\end{array}\right|z\!\right),
\end{equation}
where the function ${_{2}\Psi_1}$ on the right hand side possesses the integral representation given by the left hand side of \eqref{eq:2Psi1int}.  Under additional condition $\omega<6$  we can expand ${_{2}\Psi_1}$ from \eqref{eq:pPsiq2Psi1} by  Theorem~\ref{th:FoxWrightrho}.  However, we were unable to demonstrate that this expansion substituted  to \eqref{eq:pPsiq2Psi1} gives two convergent series and all the more that interchange of the order of summations would be legitimate.  Hence, we have

\begin{open}
Find convergent or at least asymptotic expansion for ${_{2}\Psi_1}(\rho{z})$ in the neighborhood of $z=1$ if $0<\min(\A,\B)<1/6$.
\end{open}

Finally, we mention that the principal terms of the asymptotic approximations for the Gauss hypergeometric function ${}_2F_1$ generalized in Theorems~\ref{th:FoxWrightrho} and Theorem~\ref{th:FoxWrightlogcase} were demonstrated to be monotone in \cite{ABRVV,PonVuorinen} under certain restrictions  on parameters.  More precisely, it was shown that the ratio of and/or difference between the function ${}_2F_1$ and the principal asymptotic term behaves monotonically as $z$ approaches $1$ along the real interval $(0,1)$.  These results were recently partially generalized to ${}_{p+1}F_{p}$, $p\ge2$, in \cite{WCS}.
This leads to our next open problem.
\begin{open}
Establish conditions on parameters of ${_{p}\Psi_q}(\rho{z})$  ensuring that the principal terms of the expansions \eqref{eq:pPsiq2}, \eqref{eq:pPsiqlogpos}, \eqref{eq:pPsiqlognegative} approximate this function monotonically as $z$ approaches $1$ along the real interval $(0,1)$.
\end{open}

For example, following the analogy with \cite{ABRVV}, in the zero-balanced case Corollary~\ref{cr:pPsiq-zerobal} motivates considering the monotonicity of the functions
$$
f_1(x)={_{p}\Psi_q}\left(\left.\!\!\begin{array}{c}(\a,\A)\\(\b,\B)\end{array}\right|\rho{x}\!\right)+V_0\log(1-x)
$$
and
$$
f_1(x)=\frac{x}{\log(1-x)}{_{p}\Psi_q}\left(\left.\!\!\begin{array}{c}(\a,\A)\\(\b,\B)\end{array}\right|\rho{x}\!\right).
$$


\end{document}